\documentclass[12pt,leqno]{article}
\usepackage[a4paper, margin=25.4mm]{geometry}
\usepackage{amssymb}
\usepackage{amsmath}
\usepackage{amsthm}
\usepackage{stmaryrd}
\usepackage{mathtools,mathrsfs}
\usepackage[pagebackref]{hyperref}
\hypersetup{citecolor=purple, linkcolor=blue, colorlinks=true}
\usepackage{booktabs}
\usepackage{enumerate}
\usepackage{tikz}
\usepackage{url}
\usepackage[T2A,T1]{fontenc}
\usepackage{authblk}  

\newtheorem{theorem}{Theorem}[section]

\newtheorem{lemma}[theorem]{Lemma}
\newtheorem{corollary}[theorem]{Corollary}
\theoremstyle{definition}
\newtheorem{definition}[theorem]{Definition}

\newtheorem{remark}[theorem]{Remark}

\renewcommand{\le}{\leqslant}
\renewcommand{\ge}{\geqslant}

\newcommand{\coloneq}{\vcentcolon=}      

\newcommand{\Char}{\textup{char}}

\newcommand{\F}{{\mathbb F}}

\newcommand{\Gal}{\textup{Gal}}

\newcommand{\Ord}{\textup{ord}}

\newcommand{\Tr}{\textup{Tr}}

\newcommand\kplus[1]{\kern-#1pt+\kern-#1pt}

\definecolor{mygreen}{rgb}{0.55, 0.71, 0.0}
\definecolor{myyellow}{rgb}{1.0, 0.88, 0.21}
\definecolor{mywhite}{gray}{1.00}
\definecolor{mygrey}{gray}{0.92} 

\author[1]{S.P. Glasby}
\affil[1]{\small Center for the Mathematics of Symmetry and Computation,
  University of\newline Western Australia, Perth 6009, Australia\quad
  \href{mailto:Stephen.Glasby@uwa.edu.au}{Stephen.Glasby@uwa.edu.au}}

\title{Hilbert's Theorem 90, periodicity, and roots\\ of Artin-Schreier polynomials}

\date{\today}

\begin{document}
\maketitle

\begin{abstract}
  Let $E/F$ be a cyclic field extension of degree~$n$, and let
  $\sigma$ generate the group $\Gal(E/F)$. If
  $\Tr^E_F(y)=\sum_{i=0}^{n-1}\sigma^i y=0$, then the additive form of
  Hilbert's Theorem~90 asserts that $y=\sigma x-x$ for some $x\in E$.
  When $E$ has characteristic $p>0$ we prove that $x$ gives rise to
  a periodic sequence $x_0,x_1,\dots$ which has period $pn_p$, where
  $n_p$ is the largest $p$-power that divides~$n$.
  We also show, if $y$ lies in the finite field $\F_{p^n}$, then the roots
  of a reducible Artin-Schreier polynomial $t^p-t-y$ have the form $x+u$
  where $u\in\F_p$ and 
  $x=\sum_{i=0}^{n-1}\sum_{j=0}^{i-1}z^{p^j}y^{p^i}$ for some
  $z\in\F_{p^e}$ with $e=n_p$. Furthermore, the sequence
  $\left(\sum_{j=0}^{i-1}z^{p^j}\right)_{i\ge0}$ is periodic with period $pe$.
   \vskip2mm\noindent 
  {\bf Keywords:} Hilbert's Theorem 90, Artin-Schreier polynomials, periodic sequences
  \vskip2mm\noindent
  {\bf 2020 Mathematics Subject Classification:} 11T06, 12F10, 11T55
\end{abstract}


\section{Introduction}\label{S:I}

Galois was initially concerned to determine when the roots of polynomial
$f(x)$ over a field $F$ of characteristic zero could be expressed as a
formula involving the symbols $+, -, \times, /$ and $n$th roots for $n\ge2$.
If $F$ has prime characteristic, then the formulas
were adapted to included roots of Artin-Schreier polynomials. This note
explores different formulas for the roots, replacing roots of Artin-Schreier
polynomials with polynomial expressions.

Given a quadratic equation $ax^2+bx+c=0$ over a finite
field $\F_{2^n}$ of order~$2^n$, Chen~\cite{Chen} gives `formulas' for
the roots that involve $a,b,c$ and some element of trace 1. 
We give a unified treatment of these formulas, which in prime
characteristic $p$ depends on $n_p$,
the largest $p$-power divisor of $n$. This is related to
the additive version of Hilbert's Theorem~90 and to sequences of
period~$pn_p$. We say that a sequence $(x_i)_{i\ge0}$ has \emph{period}~$e$
if $x_i=x_j$ whenever $i\equiv j\pmod{e}$, and $e$ is minimal.
We prove that a root of a reducible polynomial
$x^p-x-y\in\F_{p^n}[x]$ has the form $x=\sum_{i=0}^{n-1}x_iy^{p^i}$,
$x_i\in\F_{p^n}$, where the sequence $(x_i)_{i\ge0}$ has period $pn_p$.
In the case that $p=2$ and $n$ is odd, we may take $x_i=i\in\{0,1\}$, and
if $p=2$ and $n$ is twice odd, we may take $x_i=0,\omega,1,\omega^2$ whenever
$i\equiv0,1,2,3\pmod 4$: here $\omega$ is a primitive cube root of~1.
The periodic sequence $(x_i)_{i\ge0}$ is described in detail, 
for larger values of the 2-part $n_2$ of $n$, in Lemma~\ref{L:F2}.

This note proves the following two main results which we contextualize below.
Given a finite Galois extension $E/F$ with group $\Gamma$, the trace map
$E\to F$ is $\Tr^E_F(y)=\sum_{\gamma\in\Gamma}\gamma y$.

\begin{theorem}\label{T:H90}
  Let $E/F$ be a cyclic extension of degree~$n$. Suppose
  $\Gal(E/F)=\langle\sigma\rangle$ and $y,z\in E$ satisfy
  $\Tr^E_F(y)=0$, $\Tr^E_F(z)=1$, then $y=\sigma x-x$ where
  $x=\sum_{i=0}^{n-1}\kern-1.1pt\left(\sum_{j=0}^{i-1}\sigma^j z\right)\kern-1.1pt\sigma^i y$.
\end{theorem}

The element $z$ in Theorem~\ref{T:H90} satisfies additional properties
in prime characteristic.

\begin{theorem}\label{T:Hilbert90}
  Let $E/F$ be a cyclic extension of degree~$n$ where
  $\Gal(E/F)=\langle\sigma\rangle$ and $\Char(F)=p$. Let $n$ have
  $p$-part $n_p$, so $n_p\mid n$ and $pn_p\nmid n$. If $y\in E$ and $\Tr^E_F(y)=0$,~then
  \begin{enumerate}[{\rm (a)}]
  \item $F(z)$ contains the fixed field of $\sigma^{n_p}$, and
    $z$ may be chosen so that $|F(z):F|=n_p$;
  \item the sequence $(\sum_{j=0}^{i-1}\sigma^j z)_{i\ge0}$ is periodic with
    period~$pe$ where $e=|F(z):F|$ and $n_p\mid e$.
  \end{enumerate}
\end{theorem}

Polynomials of the form $t^p-t-y$ where $y\in E$ were studied by Artin
and Schreier in the context of generalizing Galois theory to characteristic~$p$:
they are analogues of $p$th roots in characteristic~0.
Our theorems relate to roots of these Artin-Schreier
polynomials. Suppose $\Char(E)=p>0$ and the
$p$th power map $\sigma\colon x\mapsto x^p$ of $E$
generates $\Gal(E/F)$, e.g. when $E=\F_{p^n}$ and $F=\F_p$. If $\Tr^E_F(y)=0$
and $\Tr^E_F(z)=1$, then $R(y,z)=\sum_{i=0}^{n-1}\sum_{j=0}^{i-1}z^{p^j}y^{p^i}$
is a root of the polynomial $t^p-t-y$. For example, if $n_p=1$, we may
take $z=\frac{1}{n}$ and $x_i\coloneq\sum_{j=0}^{i-1}z^{p^j}=\frac{i}{n}$ defines a
sequence of period~$p$, and $x=\sum_{i=0}^{n-1}\frac{i}{n}y^{p^i}$ is a
root of $t^p-t-y$. Because different choices of $n_p$ involve
different choices for $z$ with $\Tr^E_F(z)=1$, this suggests that there
are infinitely many (polynomial) formulas for the roots of polynomials
$t^p-t-y\in E[t]$. If $E=\F_{p^n}$, $F=\F_p$ and $e=n_p$, then finding $z$
in $\F_{p^e}$ with $\Tr^E_F(z)=1$ is easy in practise: with
probability $1-\frac{1}{p}$ a uniformly
random element $\zeta\in\F_{p^e}$ has $\tau\coloneq\Tr^E_F(\zeta)\ne0$,
so $z=\tau^{-1}\zeta$ has $\Tr^E_F(z)=1$.

There are two versions of Hilbert's Theorem~90: the additive form
and the multiplicative form: we abbreviate these AH90 and MH90, respectively.
Theorem~\ref{T:H90} is different from usual version of AH90 in that the
formula for $x$ is very slightly different. This difference leads
to Theorem~\ref{T:Hilbert90}. We remark first that MH90 is so numbered because
it appears as Satz~90 in~\cite{Hilbert} on~p.~272, and Hilbert never
actually proved AH90. The standard proof of AH90
is given in Lang~\cite[p.\,215]{Lang} and it follows Emil Artin's ideas.
In Theorem~\ref{T:H90}, $x$ equals
$R(y,z)\coloneq\sum_{i=0}^{n-1}\sum_{j=0}^{i-1}\sigma^j z\sigma^i y$
where $z\in E$ satisfies $\Tr^E_F(z)=1$ and $y=\sigma x-x$ holds.
By contrast, the standard proof
of AH90 uses a different element, namely $x'=R(z,y)$ (with $y,z$ swapped) and $y=x'-\sigma x'$ holds.
The proofs of AH90 I have found, assume $\Tr^E_F(z)\ne0$ (but scaling to get $\Tr^E_F(z)=1$ is easy) and show that
\begin{align*}
  x'&=\frac{1}{\Tr^E_F(z)}\left(y\sigma(z)+(y+\sigma(y))\sigma^2(z)+\cdots+(y+\sigma(y)+\cdots+\sigma^{n-2}(y))\sigma^{n-1}(z)\right)\\
  &=\frac{1}{\Tr^E_F(z)}\sum_{i=0}^{n-1}\sum_{j=0}^{i-1}\sigma^j y\sigma^i z=\frac{1}{\Tr^E_F(z)}R(z,y)
\end{align*}
satisfies $y=x'-\sigma x'$. The point here is that writing $x$ in the
form $\sum_{i=0}^{n-1}x_i\sigma^iy$ for some $x_i\in E$, suggests studying
the sequence $(x_i)_{i\ge0}$, whereas the form of $x'$ does not.

We remark that the expression $x'=R(z,y)$ is very natural in the
context of Galois cohomology. Indeed, this form suggests itself via
Emmy Noether's equations, and these equations lead to analogous proofs
of AH90 and MH90 as we explain.
Emmy Noether generalized the \emph{method} used to prove AH90 and MH90 by
introducing 1-cocyles, to handle a noncyclic (finite) Galois group
$\Gamma\coloneq\Gal(E/F)$, see~\cite[Theorems~4.29 and 4.30]{Jacobson}.
Her result~\cite[Theorem~4.29]{Jacobson} says: given a
function $y\colon \Gamma\to E\setminus\{0\}$ satisfying
$y_{\beta\gamma}=y_\beta\beta(y_\gamma)$ for $\beta,\gamma\in \Gamma$,
then for some $z\in E$, $x=\sum_{\beta\in\Gamma} y_\beta\beta(z)$ satisfies
$y_\gamma=x\gamma(x)^{-1}$ for all $\gamma\in \Gamma$.
Her result~\cite[Theorem~4.30]{Jacobson} says: given a
function $y\colon \Gamma\to E$ satisfying
$y_{\beta\gamma}=y_\beta+\beta(y_\gamma)$ for $\beta,\gamma\in \Gamma$,
then for some $z\in E$, $x=\Tr^E_F(z)^{-1}\sum_{\beta\in\Gamma} y_\beta\beta(z)$
satisfies $y_\gamma=x-\gamma(x)$ for all $\gamma\in \Gamma$. When
$\Gamma=\langle\sigma\rangle$ is cyclic of order~$n$, 
the 1-cocycle conditions imply that
$y_{\sigma^i}=y(\sigma y)\cdots(\sigma^{i-1} y)$ for $0\le i<n$
where $y_1=y$, and similarly in the additive case
$y_{\sigma^i}=y+(\sigma y)+\cdots+(\sigma^{i-1} y)$
for $0\le i<n$ where $y_1=y$. This gives rise to the
multiplicative formula
$x=\sum_{i=0}^{n-1}\left(\prod_{j=0}^{i-1}\sigma^jy\right)\sigma^iz$ and
when $\Tr^E_F(z)=1$, the additive formula
$x=R(z,y)\coloneq\sum_{i=0}^{n-1}\left(\sum_{j=0}^{i-1}\sigma^jy\right)\sigma^iz$, respectively.
Theorem~\ref{T:H90} uses the related expression 
$R(y,z)=\sum_{i=0}^{n-1}\left(\sum_{j=0}^{i-1}\sigma^jz\right)\sigma^iy$.
The asymmetry between $R(y,z)$ and $R(z,y)$ is measured by the
equation $R(y,z)+R(z,y)=-\sum_{i=0}^{n-1}\sigma^i(yz)$, see
Lemma~\ref{L:Ryz} below.

\section{Background}\label{S:Background}

This work fits into a very classical setting, namely generalizing Galois
theory to characteristic~$p$.
Theorem~\ref{T:Artin-Schreier} below follows from
the the additive form (AH90) of Hilbert's Theorem~90, it is proved in~\cite[Theorem~4.33]{Jacobson}
and \cite[VIII,\,Theorem~11]{Lang}.

\begin{theorem}[Artin-Schreier]\label{T:Artin-Schreier}
  Let $E$ be field of prime characteristic $p$.
  \begin{enumerate}[{\rm (a)}]
  \item Let $E/F$ be a cyclic extension of degree~$p$. Then $E=F(x)$
    where $x^p-x\in F\setminus\{0\}$.
  \item Given $y\in E$, the polynomial $t^p-t-y$ is either irreducible over $E$,
    or it has $p$ distinct roots $x+u$ where $x\in E$ and $u\in\F_p$.
  \end{enumerate}
\end{theorem}

The trace map $\Tr_F^E(x)=\sum_{\gamma\in\Gal(E/F)}\gamma(x)$ of a
finite Galois extension $E/F$, is an $F$-linear map $E\to F$ which is nonzero
by Artin's Lemma~\cite[VIII,\,Theorem~7]{Lang}. In certain cases, the trace map
can be used to distinguish whether or not $t^p-t-y$ is irreducible.

\begin{lemma}\label{L:Artin-Schreier}
  Let $E/F$ be a cyclic Galois extension where $|E:F|=\Char(F)=p$. Suppose that
  $\Gal(E/F)=\langle\sigma\rangle$ where $\sigma\colon E\to E$ is the
  $p$th power map $\sigma(x)=x^p$. For $y\in E$, the polynomial
  $t^p-t-y$ is
  reducible over~$E$ if and only if $\Tr_F^E(y)=\sum_{i=0}^{p-1}y^{p^i}=0$.
  Furthermore, $E=F(z)$ for some $z\in E$ where
  $z^p-z^{p-1}+1=0$ and $\Tr^E_F(z)=1$.
\end{lemma}

\begin{proof}
  If $\Tr_F^E(y)=0$, then  $y=\sigma(z)-z=z^p-z$
  for some $z\in E$ by AH90. Hence $t^p-t-y$ is reducible over~$E$.
  Conversely, suppose if $t^p-t-y$ is reducible over~$E$, then $t^p-t-y$
  splits over $E$ by Theorem~\ref{T:Artin-Schreier}(b). Let $z\in E$ be
  a root. Using $\sigma z=z^p$ shows
  \[
  0=\Tr^E_F(0)=\Tr^E_F(z^p-z-y)=\Tr^E_F(\sigma z)-\Tr^E_F(z)-\Tr^E_F(y)=-\Tr^E_F(y).
  \]
  Finally, set $y=-1$. Then $\Tr_F^E(-1)=p(-1)=0$ and by the above argument
  $t^p-t+1=0$ is reducible over $E$. Hence some $\zeta\in E$ satisfies
  $t^p-t+1=0$ by Theorem~\ref{T:Artin-Schreier}(b). However,
  $\sigma\zeta=\zeta^p=\zeta-1$ implies $\sigma^i\zeta=\zeta-i$ so
  $\zeta$ has $p$ distinct conjugates in~$E$ and $E=F(\zeta)$. Therefore
  $E=F(z)$ where $z=\zeta^{-1}$ is a root of the (reverse) polynomial
  $t^p(t^{-p}-t^{-1}+1)=t^p-t^{p-1}+1$. Since the conjugates
  $z,\sigma z,\dots,\sigma^{p-1}z$ are pairwise distinct,
  $t^p-t^{p-1}+1$ is irreducible over~$F$. 
  The coefficient of $t^{p-1}$ shows that $\Tr^E_F(z)=\sum_{i=0}^{p-1}\sigma^iz=1$.
\end{proof}

We note that only part of
Theorem~\ref{T:Artin-Schreier}(b) generalizes. If $\F_q\le F\le E$,
and $t^q-t-y$ has a root $x\in E$, then
$t^q-t-y$ factors as $\prod_{a\in\F_q}(t-(x+a))$.
However, if $t^q-t-y$ has no root in $E$, we can
not conclude that $t^q-t-y$ is irreducible over $E$. For example, if $q=4$ then
$t^4+t+1$ factors as $(t^2+t+\omega)(t^2+t+\omega^2)$ over $\F_4=\F_2(\omega)$.

A polynomial over a field
of characteristic $p$ is solvable by radicals if its roots can be expressed
as a finite combination of $+,-,\times,\div,\sqrt[r]{\ \;}$ where $r$ is a prime
different to $p$, \emph{and} roots of Artin-Schreier
polynomials~\cite[p.\,217]{Lang}.
Note that \emph{every} irreducible polynomial of degree~$n$ over a
finite field $\F_q$,
is trivially solvable by radicals because each nonzero element of $\F_{q^n}$ is a
primitive $m$th root of unity for some $m\mid (q^n-1)$. It makes sense,
therefore, to allow non-radical expressions for the roots of polynomials over
finite fields. One such example is 
$R(y,z)=\sum_{i=0}^{n-1}\sum_{j=0}^{i-1}z^{p^j}y^{p^i}$, which involves the
`auxiliary' element~$z$.

\section{Periodicity and Hilbert's additive Theorem~90}\label{S:H90}

  Let $E/F$ be a Galois extension with group $\Gamma=\Gal(E/F)$.
  Thus $\Tr^E_F(\sigma x)=\Tr^E_F(x)$ for $\sigma\in\Gamma$. Hence if $y\in E$ has the
  form $\sigma x-x$ for some $x\in E$ and $\sigma\in\Gamma$,
  then $\Tr^E_F(y)=0$. If $\Gamma$ is \emph{finite} and \emph{cyclic}, then
  the converse holds by AH90.
  In this section $E/F$ is a cyclic extension of degree~$n$,
  and $\Gal(E/F)=\langle\sigma\rangle$.
  Thus $\sigma^n=1$ and $\Tr_F^E(y)=\sum_{i=0}^{n-1}\sigma^i y$.

\begin{proof}[Proof of Theorem~$\ref{T:H90}$]
  Suppose that $x$ has the form $\sum_{i=0}^{n-1}x_i\sigma^iy$ where
  $x_0,\dots,x_{n-1}\in E$.  We aim to solve for the coefficients~$x_i$.
  Observe that
  \[
  \sigma x-x=-x_0y+\sum_{i=1}^{n-1}(\sigma x_{i-1}-x_i)\sigma^i y
  +\sigma x_{n-1}\sigma^ny.
  \]
  Assume for the moment that there exists
  an element $z\in E$ satisfying
  \begin{equation}\label{E:z}
    \sigma x_{i-1}-x_i=-z
    \quad\textup{for $1\le i<n$, and}\quad \sigma x_{n-1}-x_0=1-z.
  \end{equation}
  Substituting these equations and using $\sigma^n y=y$ gives
  \[
  \sigma x-x=-x_0y-z\sum_{i=1}^{n-1}\sigma^iy+(1-z+x_0)y
  =y-z\sum_{i=0}^{n-1}\sigma^iy=y
  \]
  as desired. Next, assume that $z$ and $x_0$ are given. The equations
  $\sigma x_{i-1}-x_i=-z$ show that
  $x_i=\sigma^i x_0+\sigma^{i-1}z+\cdots+\sigma z+z$ holds for $1\le i<n$.
  Rewrite the $n$ equations~\eqref{E:z} as
  $\sigma^{n-i+1} x_{i-1}-\sigma^{n-i}x_i=-\sigma^{n-i}z$
  for $1\le i<n$ and $\sigma x_{n-1}-x_0=1-z$. Adding these equations gives
  $0=\sigma^n x_0-x_0=1-\Tr_F^E(z)$. Hence $\Tr_F^E(z)=1$ must hold.

  To see that these necessary equations are sufficient, choose $z\in E$ so
  that $\Tr_F^E(z)=1$ which is possible as the trace map $\Tr_F^E\colon E\to F$
  is surjective. Next, define $x_i=\sum_{j=0}^{i-1}\sigma^j z$ 
  for $0\le i<n$ where the empty sum $x_0=\sum_{j=1}^0\sigma^j z$ is $0$
  by convention. Then the $n$ equations in~\eqref{E:z} hold, and $\sigma x-x$
  equals $y$, as claimed.
\end{proof}

\begin{lemma}\label{L:Ryz}
  Let $E/F$ be a cyclic Galois extension of degree~$n$ where
  $\Gal(E/F)=\langle\sigma\rangle$. If $y,z\in E$ have traces $\Tr^E_F(y)=0$
  and $\Tr^E_F(z)=1$, then
  $R(y,z)=\sum_{i=0}^{n-1}\left(\sum_{j=0}^{i-1}\sigma^j z\right)\sigma^i y$
  satisfies
  $R(y,z)=-\Tr^E_F(yz)-R(z,y)$, $\sigma R(y,z)-R(y,z)=y$ and
  $R(z,y)-\sigma R(z,y)=y$.
\end{lemma}

\begin{proof}
  The equations $\sum_{i=0}^{n-1}\sigma^i y=0$ and
  $\sum_{j=0}^{n-1}\sigma^j z=1$ imply
  \[
  R(y,z)=\sum_{i=0}^{n-1}\left(\sum_{j=0}^{i-1}\sigma^j z\right)\sigma^i y
        =\sum_{i=0}^{n-1}\left(1-\sum_{j=i}^{n-1}\sigma^j z\right)\sigma^i y
        =-\sum_{i=0}^{n-1}\left(\sum_{j=i}^{n-1}\sigma^j z\right)\sigma^i y.
  \]
  By swapping the order of summation, by swapping $i\leftrightarrow j$, and
  by using the formula $\Tr^E_F(zy)=\sum_{i=0}^{n-1}\sigma^i(zy)$, the
  previous display may be written 
  \[
  R(y,z)=-\sum_{j=0}^{n-1}\sum_{i=0}^{j}\sigma^j z\sigma^i y
        =-\sum_{i=0}^{n-1}\sum_{j=0}^{i}\sigma^i z\sigma^j y
        =-\Tr^E_F(zy)-\sum_{i=0}^{n-1}\left(\sum_{j=0}^{i-1}\sigma^j y\right)\sigma^i z.
  \]
  Therefore $R(y,z)=-\Tr^E_F(yz)-R(z,y)$. Theorem~\ref{T:Hilbert90} proves
  that $\sigma R(y,z)-R(y,z)=y$ holds. Since $\sigma\Tr^E_F(yz)=\Tr^E_F(yz)$
  it follows that $R(z,y)-\sigma R(z,y)=y$, as claimed.
\end{proof}

We prove in Theorem~\ref{T:Hilbert90} that when $\Char(F)=p>0$,
the element $x$ in Theorem~\ref{T:H90} has an
intrinsic periodicity.  This observation is used in
Corollary~\ref{C:Hilbert90} to give a polynomial formula for the
roots of (reducible) Artin-Schreier polynomials.

\begin{proof}[Proof of Theorem~$\ref{T:Hilbert90}$]
  (a)~First, we show that the fixed subfield $K$ of $\sigma^{n_p}$ contains
  an element $z$ satisfying $\Tr^E_F(z)=1$.
  Let $n_{p'}\coloneq n/n_p$ be the $p'$-part of $n$.
  Since $n_{p'}\ne0$ in $F$ and the map
  $\Tr^K_F\colon K\to F$ is surjective, there exists an element $z\in K$
  such that $\Tr^K_F(z)=(n_{p'})^{-1}$. Therefore $\Tr^E_K(z)=|E:K|z=n_{p'}z$ and hence
  $\Tr^E_F(z)=\Tr^K_F(\Tr^E_K(z))=n_{p'}\Tr^K_F(z)=1$.
  Next, suppose that $z\not\in K$ and $\Tr^E_F(z)=1$. Let $L=F(z)$.
  Then $1=\Tr^L_F(\Tr^E_L(z))=|E:L|\Tr^L_F(z)$.
  Therefore $p\nmid|E:L|$, so $n_p=|K:F|$ divides $e=|L:F|$ as claimed.
  In summary, the minimal polynomial over $F$ of an element $z$ 
  with $\Tr^E_F(z)=1$ has degree a multiple of $n_p$, and degree
  precisely $n_p$ if $z$ lies in $K\setminus K_0$ where $K_0$ is the
  unique (maximal) subfield of $K$ with $F\le K_0<K$ and $|K:K_0|=p$.

  (b)~Let $x_i=\sum_{j=0}^{i-1}\sigma^j z$ for $0\le i<n$ where $\Tr^E_F(z)=1$.
  Let $L=F(z)$ and $e=|L:F|$. Hence $x_e=\Tr^L_F(z)\ne0$ by part~(b)
  and $\sigma^j z=z$ implies $e\mid j$. Suppose that $x_k=0$ where
  $0\le k<pe$ and write $k=ie+j$ where $0\le i<p$ and $0\le j<e$.
  The equation $x_{\ell e+j}=\ell x_e+x_j$ shows
  that $x_0=x_{pe}=0$. Thus $x_k=ix_e+x_j=0$ and $x_e\in F$ gives
  \[
  x_{j+1}-z=\sigma x_j=\sigma(-ix_e)=-ix_e=x_j.
  \]
  Whence $\sigma^j z=x_{j+1}-x_j=z$. This implies that
  $e\mid j$ and hence $j=0$. However, $0=x_k=x_{ie}=ix_e$ and $x_e\ne0$ shows
  that $0=i=j=k$. Therefore $x_\ell\ne0$ holds for
  $0<\ell<pe$, and the sequence $x_0,x_1,\dots$ is periodic, with period
  precisely~$pe$.
\end{proof}

\section{Formulas for the roots of Artin-Schreier polynomials}\label{S:roots}

In this section we assume that $F$ has prime characteristic $p$. Set $q=p^f$.
Two examples with $\Gal(E/F)=\langle\sigma\rangle$ cyclic of order $n$ and
$\sigma(x)=x^q$ for $x\in E$, are the finite fields $E=\F_{q^n}$, $F=\F_q$,
and the function fields $E=\F_{q^n}(u)$, $F=\F_q(u)$.
We shall give a number of different applications of Theorem~\ref{T:Hilbert90}
depending on the value of $n_p$. Observe that an immediate corollary of
Theorem~\ref{T:H90} is a polynomial formula for a root of $t^q-t-y$.

\begin{corollary}\label{C:Hilbert90}
  Suppose that $\Char(F)=p$, $q=p^f$, and $E/F$ is a cyclic extension of
  degree~$n$ where $\Gal(E/F)=\langle\sigma\rangle$ and $\sigma(x)=x^q$
  for $x\in E$.  If $y\in E$ and  $\Tr^E_F(y)=0$ and $\Tr^E_F(z)=1$, then
  $x=\sum_{i=0}^{n-1}\sum_{j=0}^{i-1}z^{q^j}y^{q^i}$ satisfies $x^q-x-y=0$.
\end{corollary}

When $n$ is coprime to $p$, i.e. $n_p=1$, the formula 
for $x$ in Corollary~\ref{C:Hilbert90} simplifies.

\begin{corollary}\label{C:n_p=1}
  Suppose that $\F_q\le F$, $q=p^f$, and $E/F$ is a cyclic extension of
  degree~$n$ where $\Gal(E/F)=\langle\sigma\rangle$ and $\sigma(x)=x^q$
  for $x\in E$.
  If $y\in E$,  $\Tr^E_F(y)=0$ and $n_p=1$, then the factorization
  $t^q-t-y=\prod_{u\in\F_q}(t-x-u)$
  holds where $x=\sum_{i=0}^{n-1}\frac{i}{n}y^{q^i}$.
\end{corollary}

\begin{proof}
  By Theorem~\ref{T:Hilbert90}(b), we may choose $z$ to lie in $F$
  and so $1=\Tr^E_F(z)=nz$. Hence $z=\frac{1}{n}$ and
  $x=R(y,z)=\sum_{i=0}^{n-1}\frac{i}{n}y^{q^i}$ satisfies $x^q-x-y=0$
  by Corollary~\ref{C:Hilbert90}. Therefore $x+u$ is a root of $t^q-t-y$ for
  each $u\in\F_q\le F$, and the claimed factorization follows.
\end{proof}
  
\begin{remark}
  When $\Char(F)=p\nmid n$, Corollary~\ref{C:n_p=1} implies
  that all the roots of $t^q-t-y$ lie
  in the $\F_p$-subspace
  $\F_q+\{\sum_{i=0}^{n-1}\lambda_iy^{q^i}\mid\lambda_1,\dots,\lambda_{n-1}\in\F_p\}$ of $\F_{q^n}$.
\end{remark}

When $n_p=p$, i.e. $p\mid n$ and $p^2\nmid n$, the element $z$
in Corollary~\ref{C:Hilbert90} can be chosen to lie in the fixed
field of $\sigma^{p}$, and to have a specific minimal polynomial.

\begin{corollary}\label{C:n_p=p}
  Let $E/F$ be a cyclic field extension of degree $n$ where $n_p=\Char(F)=p$.
  Suppose that
  $\Gal(E/F)=\langle\sigma\rangle$ where $\sigma\colon E\to E$ is the
  $p$th power map $\sigma(x)=x^p$.
  If $y\in E$ satisfies $\Tr^E_F(y)=\sum_{i=0}^{n-1}y^{p^i}=0$, and $\zeta$
  satisfies $\zeta^p-\zeta^{p-1}+1=0$. Then $x^p-x-y=0$ where
  $x=(n/p)^{-1}\sum_{i=0}^{n-1}\sum_{j=0}^{i-1}\zeta^{p^j}y^{p^i}$
  and $\zeta^{p^{p}-1}=1$.
\end{corollary}

\begin{proof}
  It follows from $n_p=p$, that $n/p$ is coprime to $p$.
  Let $K$ be the fixed field of~$\sigma^{p}$.  Then $|K:F|=p$ and $K=F(\zeta)$
  where $\zeta$ satisfies $\zeta^p-\zeta^{p-1}+1=0$ and
  $\Tr^K_F(\zeta)=1$ by Lemma~\ref{L:Artin-Schreier}.
  Therefore $\Tr^E_F(\zeta)=n/p$
  and $\Tr^E_F(z)=1$ where $z=(n/p)^{-1}\zeta$. By Theorem~\ref{T:H90}, 
  $x=(n/p)^{-1}\sum_{i=0}^{n-1}\sum_{j=0}^{i-1}\zeta^{p^j}y^{p^i}$ satisfies
  $x^p-x-y=0$. Finally, $\sigma(\zeta)=\zeta^p$ implies
  that $\zeta=\sigma^{p}(\zeta)=\zeta^{p^{p}}$, so that $\zeta^{p^{p}-1}=1$.
\end{proof}

We shall now apply Corollary~\ref{C:Hilbert90} in the case that
$E=\F_{p^n}$ and $F=\F_p$. In this case $E/F$ is cyclic and $\Gal(E/F)=\langle\sigma\rangle$ has order~$n$
where $\sigma(x)=x^p$ is the $p$th power map.

Let $p,r$ be distinct primes. Denote by $\Ord_r(p)$ the
smallest positive integer $e$ such that $p^e\equiv1\pmod r$ and note that
$e\mid(r-1)$. The cyclotomic polynomial $\Phi_r(t)=\sum_{i=0}^{r-1}t^i$
of degree $r-1$
factors over $\F_p$ as a product of $(r-1)/e$ irreducible polynomials each of
degree $e=\Ord_r(p)$ by~\cite[Theorem~2.47]{LN}. 

Theorem~\ref{T:Hilbert90} shows that $z$ may be chosen to lie in a subfield
of order $p^f$ where $f={n_p}$. Hence $z^{p^f-1}=1$.
Lemma~\ref{L:p,r} proves that, in many cases, $z$ may be chosen to
have prime order $r\ne p$. This choice of $z$, however, may increase
the period of the sequence $(x_i)_{i\ge0}$.
  
\begin{definition}\label{D:big}
  Call a non-zero polynomial $f(t)=\sum_{i=0}^d f_it^i$ \emph{big}
  if $d=0$ or $f_df_{d-1}\ne0$, and \emph{small} if $f_d\ne0$ and $f_{d-1}=0$.
  As a small polynomial times a small polynomial is small, and a big
  polynomial times a small polynomial is big, we have the following lemma.
\end{definition}

\begin{lemma}\label{L:big}
  If a big polynomial factors, then at least one of its factors is big.
\end{lemma}

\begin{lemma}\label{L:p,r}
  Let $p, r$ be unequal primes. Let $e\coloneq\Ord_r(p)$ divide $n$. Set
  $E=\F_{p^n}$, $K=\F_{p^e}$ and $F=\F_p$. There exists a primitive
  $r$th root of unity $\zeta\in K$ 
  for which $\tau\coloneq\Tr^K_F(\zeta)\ne0$.
  If $p\nmid \frac{n}{e}$ and $y\in E$ has $\Tr^E_F(y)=0$, then
  $x=(\frac{n\tau}{e})^{-1}\sum_{i=0}^{n-1}\left(\sum_{j=0}^{i-1}\zeta^{p^j}\right)y^{p^i}$
  satisfies $x^p-x-y=0$. Finally, the sequence
  $(\sum_{j=0}^{i-1}\zeta^{p^j})_{i\ge0}$ is periodic with period~$ep$
  where $n_p\mid e$.
\end{lemma}

\begin{proof}
  First, $r$ divides $p^d-1$ if and only if $e\mid d$.
  As $e\mid n$, all primitive $r$th roots of unity lie in $K=\F_{p^e}$, and
  $F\le K\le E$ holds. The cyclotomic polynomial
  $\Phi_r(t)=\sum_{i=0}^{r-1}t^i$ is big, and
  $\Phi_r(t)$ factors over $F$ as a product of irreducible  polynomials,
  each of degree~$e$, at least one of which is big by Lemma~\ref{L:big}. Hence
  some $r$th root of  unity $\zeta\in K$  has $\tau\coloneq\Tr^K_F(\zeta)\ne0$.
  Thus $\Tr^E_K(\zeta)=\frac{n}{e}\zeta$  and
  $\Tr^E_F(\zeta)=\Tr^K_F(\Tr^E_K(\zeta))=\frac{n}{e}\Tr^K_F(\zeta)=\frac{n\tau}{e}\ne0$.
  Therefore $z=(\frac{n\tau}{e})^{-1}\zeta\in K^\times$ satisfies $\Tr^E_F(z)=1$.
  The formula for $x=R(y,z)$ follows from Theorem~\ref{T:H90} and the
  periodicity $ep$ and $n_p\mid e$ follows from  Theorem~\ref{T:Hilbert90}(b).
\end{proof}

\begin{lemma}\label{C:p=2mod3}
  Let $p$ be a prime where $p\equiv2\pmod3$. Let $E=\F_{p^n}$ and $F=\F_p$ where
  $n$ is even and let $\omega\in E$ satisfy $\omega^2+\omega+1=0$.
  If $p\nmid(n/2)$, and $y\in E$ satisfies $\Tr^E_F(y)=0$, then
  the sum
  $x=(\frac{n}{2})^{-1}\sum_{i=0}^{n-1}\left(\lfloor\frac{i}{2}\rfloor-r_i\omega\right)y^{p^i}$
  satisfies $x^p-x-y=0$ where  $r_i=i-2\lfloor i/2\rfloor$.
  Finally, the sequence $x_0, x_1,\dots$ is periodic of period~$2p$ where
  $x_i=(\frac{n}{2})^{-1}\left(\lfloor\frac{i}{2}\rfloor-r_i\omega\right)$.
\end{lemma}

\begin{proof}
  Since $n$ is even and $p^n\equiv1\pmod3$, $E$ contains a primitive cube
  root~$\omega$ of 1.
  The condition $p\equiv 2\pmod3$ is equivalent to $e=\Ord_3(p)=2$. Thus we may
  use Lemma~\ref{L:p,r} with $r=3$, $\zeta=\omega$ and $\tau=-1$.
  However, $\omega^{p^k}$ equals $\omega$
  if $k$ is even, and $\omega^2$ if $k$ is odd. If $i=0,1,2,3,4,5,6$, then
  $s_i=\sum_{k=0}^{i-1}\omega^{p^k}$ equals $0$, $\omega$, $\omega+\omega^2=-1$,
  $-1+\omega$, $-2$, $-2+\omega$, $-3$, etc. Therefore
  $s_i=-(i-r_i)/2 +r_i\omega$ where $r_i=i-2\lfloor i/2\rfloor\in\{0,1\}$.
  Consequently, the values of
  $x_i=(\frac{n}{2})^{-1}\left(\frac{i-r_i}{2} -r_i\omega\right)$
  are periodic of period~$2p$ by Theorem~\ref{T:Hilbert90}(b).
\end{proof}

The following result gives different quadratic `formulas' over the finite
field~$\F_{2^n}$ that depend on $n_2$ and the choice of
an element $z\in\F_{2^n}$ with absolute trace~1. The condition $n_2=2^k$ is
equivalent to $n\equiv 2^k\pmod{2^{k+1}}$. The case $k=0$ is handled by
Corollary~\ref{C:n_p=1}.

\begin{lemma}\label{L:F2}
  Let $y\in\F_{2^n}$ where $\sum_{i=0}^{n-1}y^{2^i}=0$. Then 
  $x=\sum_{i=0}^{n-1}x_iy^{2^i}$ is a root of $t^2+t+y$ where
  $x_i=\sum_{j=0}^{i-1}z^{2^j}$ and $z\in\F_{2^n}$ has absolute trace
  $\sum_{i=0}^{n-1}z^{2^i}=1$ and has degree $n_2$ and the sequence $(x_i)$
  has period $2n_2$. Table~$\ref{T:z}$ lists choices
  for $z$  and the minimal polynomial $m_{z/\F_2}(t)$ of $z$ over $\F_2$, as a
  function of the $2$-part $n_2$ of $n$.
  \begin{table}[ht!]
    \begin{center}\label{T:z}
    \caption{Choices for $z$ and its minimal polynomial $m_{z/\F_2}(t)$ over $\F_2$ of degree $n_2$}
    \begin{tabular}{ llllll } 
    $n_2$&$2$&$4$&$8$&$16$&$32$\\ \hline
    $z$&$\omega$&$\alpha$&$\beta$&$\gamma$&$\delta$\\
    $m_{z/\F_2}(t)$&$t^2\kplus{2}t\kplus{2}1$&$t^4\kplus{2}t^3\kplus{2}1$&$t^8\kplus{2}t^7\kplus{2}t^2\kplus{2}t\kplus{2}1$&$t^{16}\kplus{2}t^{15}\kplus{2}t^8\kplus{2}t\kplus{2}1$&$t^{32}\kplus{2}t^{31}\kplus{2}t^3\kplus{2}t\kplus{2}1$
    \end{tabular}
    \end{center}
  \end{table}
\end{lemma}

\begin{proof}
  Set $E=\F_{2^n}$, $F=\F_2$ and $e=n_2$. The last line of
  Table~\ref{T:z} lists big irreducible polynomials over $F$ whose roots
  have order $2^{n_2}-1$, see Lemma~\ref{L:big} and Remark~\ref{R:z}. The
  polynomials are big by Definition~\ref{D:big}.
  The irreducibility over $F$, and the order of~$z$, are easily checked using a
  computational algebra system. These together imply that $K\coloneq F(z)$ has
  $|K:F|=n_2$ and $\Tr^K_F(z)=1$. Hence $\Tr^E_F(z)=|E:K|=n_{2'}=1$ as desired.
  By Theorem~\ref{T:H90}, $x=R(y,z)$ is a root of $t^2+t+y$.
  In order to use our formulas, the reader will need to know
  the values of $x_i=\sum_{j=0}^{i-1}z^{2^j}$. Since $x_0=x_{2e}=0$ and
  $x_{i+2e}=x_i$, it suffices to list $x_i$ for $0\le i<2e$. However,
  since $x_e=\Tr_F^Kz=1$ and $x_{i+e}=1+x_i$, it suffices to list $x_i$
  for $0\le i<e$.
  
  (a)~Suppose that $n_2=1$.
  Set $q=p=2$ in Corollary~\ref{C:n_p=1}. Then $x^2+x+y=0$ where
  $x=\sum_{i=0}^{n-1}\frac{i}{n}y^{2^i}=\sum_{j=0}^{(n-1)/2}y^{2^{2j+1}}$. Here
  $x_i\in\{0,1\}$ satisfies $x_i\equiv i\pmod 2$.

  (b)~Suppose that $n_2=2$.
  Set $p=2$ and $z=\omega$ into Lemma~\ref{C:p=2mod3}.
  Since $n/2$ is odd the value of
  $x_i=(\frac{n}{2})^{-1}\left(\lfloor\frac{i}{2}\rfloor-r_i\omega\right)$
  is $0,\omega,1,\omega^2$ when $i\equiv 0,1,2,3\pmod4$;
  recall that $r_i=i-2\lfloor\frac{i}{2}\rfloor\in\{0,1\}$.
  Equivalent formulas for $x$ in cases~(a) and~(b) were given previously by
  Chen~\cite[Theorems~1,\,2]{Chen} and by Lahtonen~\cite{L}.

  (c) Suppose that $n_2=4$. Then
  $z=\alpha\in\F_{16}$ has order $2^4-1$ and $\Tr^{K}_F(z)=1$.
  Also the values of $x_i=\sum_{j=0}^{i-1}z^{2^i}$ are
  $0, \alpha, \alpha^{13}, \alpha^6,1,\alpha^{12}, \alpha^7,\alpha^8$
  for $i\equiv 0,1,\dots,7\pmod 8$.
  
  (d) Suppose that $n_2=8$.
  Then $z=\beta\in\F_{256}$ has order $2^8-1$ and $z$ satisfies
  $\Tr^{K}_F(z)=1$. Furthermore, the values of $x_i=\sum_{j=0}^{i-1}z^{2^i}$ are
  $0, \beta, \beta^{100}, \beta^{189}, \beta^{29}, \beta^{60}, \beta^{154}, \beta^{177}$  for $i\equiv 0,1,\dots,7\pmod {16}$. Also $x_{i+8}=1+x_i$
  for $i\equiv 0,1,\dots,7\pmod {16}$.
  
  (e) Suppose that $n_2=16$.
  Then $z=\gamma\in\F_{2^{16}}$ has order $2^{16}-1$ and $\Tr^{K}_F(z)=1$.
  In addition, $x_i=\sum_{j=0}^{i-1}z^{2^j}$ equals
  $0, \gamma, \gamma^{64409}, \gamma^{48754}, \gamma^{27742}, \gamma^{48469}, \gamma^{1146}, \gamma^{22404}, \gamma^{64313}, \gamma^{47682}, $  $\gamma^{63219}, \gamma^{45929}, \gamma^{55680},\gamma^{46875}, \gamma^{7495}, \gamma^{32204}$ and $x_{i+16}=1+x_i$ for $i\equiv 0,1,\dots,15\pmod {32}$.

  (f) Suppose that $n_2=32$.
  Then $z=\delta\in\F_{2^{32}}$ has order $2^{32}-1$ and $\Tr^{K}_F(z)=1$.
  Writing $x_i=\sum_{j=0}^{i-1}z^{2^j}$ as a power of $\delta$ for
  $i\equiv 0,1,\dots,15\pmod {32}$, involves solving the discrete logarithm
  problem.
  This was done in the previous cases as writing a power of a primitive element
  was shorter than writing polynomials in the primitive element.
\end{proof}

\begin{remark}\label{R:z}
  It is an open question whether big primitive polynomials of
  degree $e=n_2=2^k$ in the last line of Table~\ref{T:z}
  exist for arbitrarily large choices of $k$. This turns out to be true
  if $2^{e}-1$ is square-free because by Lemma~\ref{L:tensor} the cyclotomic
  polynomial $\Phi_{2^e-1}(t)$ is big. This follows because $\Phi_{2^e-1}(t)$
  is a tensor product~\cite{G} 
  $\Phi_{p_1}(t)\otimes\cdots\otimes\Phi_{p_r}(t)$ of big polynomials
  where $2^{e}-1=p_1\cdots p_r$ and the $p_i$
  are distinct primes. However, it is not known whether
  $2^{2^k}-1=\prod_{i=0}^{k-1}(2^{2^i}+1)$ is square-free for all $k\ge1$.
  (The Fermat numbers $2^{2^i}+1$ are pairwise coprime, but it not known whether
  Fermat numbers are square-free.)
  The requirement that the polynomial be primitive (i.e. a root $z$
  has order $2^{e}-1$) is not necessary, we assumed primitivity only
  so that the non-zero elements $x_i$ could be written as a power of~$z$.
  What surprised the author is that a root of the Schreier polynomial $t^2+t+y$
  in $\F_{2^n}$ can always be described using only ${\rm O}(e^2)$ bits of~information.
\end{remark}

\begin{lemma}\label{L:tensor}
  The tensor product of big nonzero polynomials is big, see Definition~$\ref{D:big}$.
\end{lemma}

\begin{proof}
  Let $a(t)$ and $b(t)$ be big nonzero polynomials of respective degrees $m$ and $n$.
  The result is clear if $mn=0$ by~\cite[p.\,307]{G}.
  Suppose that $mn>0$ and $a(t)=\prod_{i=1}^m(t-\alpha_i)$ and
  $b(t)=\prod_{j=1}^n(t-\beta_j)$ are \emph{monic} with roots
  $\alpha_i, \beta_j$ in an extension field.
  By Definition~\ref{D:big},  $a(t)=t^m+a_{m-1}t^{m-1}+\cdots$
  and $b(t)=t^n+b_{n-1}t^{n-1}+\cdots$ where
  $a_{m-1}=-\sum_{i=1}^m\alpha_i\ne0$ and $b_{n-1}=-\sum_{j=1}^n\beta_j\ne0$.
  However,
  $(a\otimes b)(t)=\prod_{i,j}(t-\alpha_i\beta_j)$ by~\cite[p.\,307]{G}, and
  hence  $(a\otimes b)(t)=t^{mn}-(\sum_{i,j}\alpha_i\beta_j)t^{mn-1}+\cdots$.
  Therefore $(a\otimes b)(t)$ is big since
  $\sum_{i,j}\alpha_i\beta_j=a_{m-1}b_{n-1}\ne0$.
  More generally, if $a(t)$ and $b(t)$ have
  leading terms $a_m$ and $b_n$ respectively, then $(a\otimes b)(t)$ has leading
  term $a_m^nb_n^m\ne0$ by~\cite[p.\,307]{G}, and hence is big.
\end{proof}

\end{document}